\documentclass[12pt]{article}
\usepackage{fancyhdr}
\usepackage[centertags]{amsmath}
\usepackage{amsfonts}
\usepackage{graphics}
\usepackage{graphicx}
\usepackage{wrapfig}
\usepackage{amssymb}
\usepackage{amsthm}
\usepackage{newlfont}
\usepackage{url}
\usepackage{epsfig}
\usepackage{ifthen}
\usepackage{ifpdf}
\usepackage{rotating}
\usepackage{caption}
\usepackage{amsmath}
 \usepackage{subfloat}
\input amssym.def
\input amssym
\input cyracc.def

\newtheorem{lemma}{Lemma}[section]

\newtheorem{te}{Theorem}[section]

\newcommand{\be}{\begin{equation}}
\newcommand{\ee}{\end{equation}}
\newcommand{\ba}{\begin{array}}
\newcommand{\ea}{\end{array}}
\newcommand{\bee}{\begin{eqnarray*}}
\newcommand{\eee}{\end{eqnarray*}}
\newcommand{\bea}{\begin{eqnarray}}
\newcommand{\eea}{\end{eqnarray}}

\newcommand{\wmin}{w_{min}}
\newcommand{\wmax}{w_{max}}

\newcommand{\1}{\textbf{1}}

\newcommand{\jc}{\cal{J}}

\newcommand{\tx}{\tilde{x}}
\newcommand{\bx}{\bar{x}}
\newcommand{\tu}{\tilde{u}}
\newcommand{\bu}{\bar{u}}

\newcommand{\eb}{\bar{e}}
\newcommand{\si}{\Sigma_k^{-1}}
\newcommand{\dmin}{{d_{\min}}}
\newcommand{\dmax}{{d_{\max}}}
\newcommand{\ddif}{\Delta}

\begin{document}

\title{
Exact Spectral - Like Gradient Method  for Distributed Optimization
}

\author{Du\v{s}an Jakoveti\'c \footnote{Department of Mathematics and Informatics, Faculty of Sciences, University of Novi
  Sad, Trg Dositeja Obradovi\'ca 4, 21000 Novi Sad, Serbia. e-mail:
   \{djakovet@uns.ac.rs, natasak@uns.ac.rs, natasa.krklec@dmi.uns.ac.rs\}. Research supported by the Serbian Ministry of Education,
Science, and Technological Development, Grant no. 174030. 
 The work was also supported in part by the European Union (EU) Horizon
2020 project I-BiDaaS, project number 780787.
}
 \and
Nata\v sa  Kreji\'c \footnotemark[1]
\and
Nata\v{s}a Krklec Jerinki\'c \footnotemark[1]}

\maketitle

\vspace{-1cm}

 \begin{abstract}
Since the initial proposal in the late 80s, spectral gradient methods continue to receive significant attention, especially due to their excellent numerical performance on various large scale applications. However, to date, they have not been sufficiently explored in the context of distributed optimization. In this paper, we consider unconstrained distributed optimization problems where $n$ nodes  constitute an arbitrary connected network and collaboratively minimize the sum of their local convex cost functions. In this setting, 
building from existing exact distributed gradient methods, we propose a novel exact distributed gradient method wherein nodes' step-sizes are designed according to the novel rules akin to those in spectral gradient methods. We refer to the proposed method as Distributed Spectral Gradient method (DSG).
  The method exhibits R-linear convergence under standard assumptions
for the nodes' local costs and safeguarding on the algorithm step-sizes. We illustrate the method's performance through simulation examples.

\vspace{3mm}

[2mm]


\end{abstract}

\section{Introduction}

We consider a connected network with $ n $ nodes, each of which has access to a local cost function $ f_i: \mathbb{R}^d \rightarrow \mathbb{R}, \; i=1,\ldots,n. $ The objective for all nodes is to minimize the aggregate cost function  $ f: \mathbb{R}^d \rightarrow \mathbb{R} $, defined by

\be \label{objective}
f(y) = \sum_{i=1}^{n} f_i(y).
\ee

Problems of this form attract a lot of scientific interest as they arise in many emerging applications like
 distributed inference in sensor networks \cite{RibeiroADMM1, SoummyaEst, SayedConf,
SayedEstimation}, distributed control, \cite{JoaoMotaMPC},
distributed learning, e.g., \cite{BoydADMM}, etc.
 To solve this and related problems several distributed first
 order methods, e.g., \cite{nedic_T-AC,SayedEstimation,arxivVersion}, and
 second order methods, e.g., \cite{ribeiroNNpart1,DQM, DQNSIAM},
 have been proposed.
 The methods of this type converge to
 an approximate solution of problem~\eqref{objective}
 if a constant (non-diminishing) step size is used; they can be interpreted through a penalty-like
 reformulation of~\eqref{objective}; see \cite{cdc-submitted,ribeiroNNpart1}
  for details. Convergence to
  an exact solution can be achieved by using
  diminishing step-sizes, but this
  comes at a price of slower convergence.

More recently,
   \emph{exact} distributed
   first order methods, e.g., \cite{extra,harnessing,dusan},
   and second order methods~\cite{ESOM} have been proposed,
   that converge to the exact solution
   under \emph{constant} step sizes. The method in \cite{extra} uses two different weight matrices,
   differently from the standard distributed gradient method
   that utilizes a single weight matrix. The methods in~\cite{harnessing, angelia, angelia2}
   implement tracking of the network-wide average gradient and
   correct the dynamics of the standard distributed method \cite{nedic_T-AC}
   by replacing the nodes' local gradients
   with the tracked global average gradient estimates.
   A unification of the above methods
   and some further improvements are presented in~\cite{dusan}.
   An exact distributed second order method
   has been developed in~\cite{ESOM}. We refer to~\cite{dusan}
   for a detailed review of other works on exact distributed methods.

\emph{Spectral gradient methods} are a popular class of methods in centralized optimization due to their simplicity and efficiency. The class originated with the proposal of the Barzilei-Borwein method~\cite{bb}
and its analysis therein for convex quadratic functions,
while the method has been subsequently extended to
more general optimization problems, both unconstrained and constrained, \cite{raydan1993,raydan1997, bmr4}.
 Spectral gradient methods
 can be viewed as a mean to incorporate
 second-order information in a computationally efficient manner into
  gradient descent methods.
  In practice, they achieve significantly
  faster convergence with respect to
  standard gradient methods while the additional computational
  overhead per iteration is very small.
  Roughly speaking, the main idea
  of  spectral gradient methods is to approximate the Hessian at each iteration with a scalar matrix
   (the leading scalar of the matrix is called the spectral coefficient) that approximately fits the secant equation.
  Calculating the spectral gradient's scalar matrix
  is much cheaper than evaluation of the Newton direction while the convergence speed
  is usually much better than that of the gradient method.
  Spectral methods are characterized by a non-monotone behaviour which makes them suitable for combination with non-monotone line search methods, \cite{raydan1997}. It was demonstrated in \cite{raydan1997} that the spectral gradient method can be more efficient than the conjugate gradient method for
  certain classes of optimization problems.
  The R-linear convergence rate was established in~\cite{dl}, while
  extensions to constrained optimization in the form of Spectral Projected Gradient (SPG) methods are developed in \cite{bmr1,bmr2,bmr3}. A vast number of applications is available in the literature, and a comprehensive overview is presented in \cite{bmr4}.

The principal aim of this paper is to provide a generalization of spectral
gradient methods to distributed optimization 
and give preliminary numerical tests of its efficiency.
Extension of spectral gradient methods
to a distributed setting is a highly nontrivial task.
We develop an exact method (converging to the exact solution) that we refer to 
as Distributed Spectral
Gradient method (DSG). The method
utilizes step-sizes that are akin to
those of centralized
spectral methods. The
spectral-like step-sizes
 are embedded into
 the exact distributed
 first order method in~\cite{harnessing}; see also~\cite{angelia,angelia2}. 
  We utilize the primal-dual interpretation
 of the method in~\cite{harnessing} -- as provided in~\cite{dusan} (see also~\cite{angelia2}) -- 
  and the corresponding
 form of the error recursion equation.
 An analogy
 with the error recursion
 of the conventional  spectral method stated in~\cite{raydan1993} is exploited  to
 define the time-varying,
 node dependent, algorithm step-sizes.
 This analogy also allows
 for an intuitive interpretation
 of the proposed method.
 The DSG method converges R-linearly to the 
 exact solution 
 under appropriate conditions
 on the cost functions and under 
 safeguarding on the step-sizes. Initial numerical results 
 show that DSG exhibits a significantly improved convergence speed 
 with respect to its ``baseline'' method \cite{harnessing}. 
 That is, incorporating spectral-like step-sizes 
 continues to bring improvements in distributed optimization as well. 

The paper is organized as follows. Some preliminary considerations and assumptions are presented in Section~2. The proposed distributed spectral method is introduced in Section~3, while the convergence theory is developed in Section~4. Initial numerical tests are presented in Section~5, and some conclusions are drawn in Section~6.

\section{Model and preliminaries}

The network and optimization models
that we assume are described in Subsection 2.1. The proposed method is based on  the distributed gradient method developed in \cite{harnessing} and
the centralized spectral gradient method \cite{raydan1993} which are briefly reviewed in Subsection 2.2 and 2.3. The convergence analysis is based on the Small Gain Theorem which is stated in  Subsection~{2.4}.

\subsection{Optimization and network models}
We impose a set of standard assumptions on the functions  $ f_i $ in (\ref{objective}) and on the underlying network.

\noindent {\bf Assumption A1.} Assume that each function $f_{i} : \mathbb{R}^d \rightarrow \mathbb{R}, \; i=1,\ldots,n$
is twice continuously differentiable, and there exist constants $0<\mu\leq L<\infty$, such that, for every $x \in \mathbb{R}^d$, there holds
$$\mu I \preceq \nabla^{2} f_{i}(x) \preceq L I.$$
 Here, $I$ denotes the $d \times d$ identity matrix, and notation~$P  \preceq Q$ means that the matrix $Q-P$ is positive semi-definite.

Assumption~{A1} implies that each function $ f_i , \; i=1,\ldots,n $, is strongly convex with modulus $ \mu > 0, $
i.e., there holds
\be \label{A11}
f_{i}(z) \geq f_{i}(y) + \nabla f_{i}(y)^T (z-y) + \frac{\mu}{2} \|z-y\|^2, \; y,z \in \mathbb{R}^d.
\ee
Also, the gradients of the $f_i$'s are Lipschitz continuous with constant $ L, $
\be
\label{A12}
\|\nabla f_i(y) - \nabla f_i(z)\| \leq L \|y-z\|, \; y,z \in \mathbb{R}^d, \; i=1,\ldots, n.
\ee
Under Assumption~A1, problem \eqref{objective} is
solvable and has a unique solution,
denoted by~${y^*}$.
For future reference,  let us introduce the function $F:\,{\mathbb R}^{nd} \rightarrow \mathbb R$,
defined by:
\be \label{dobjective}
F(x) = \sum_{i=1}^{n} f_i(x_i),
\ee
where $x\in {\mathbb R}^{nd}$
consists of $n$ blocks $x_i \in {\mathbb R}^d$,
i.e., $x=((x_1)^T,...,(x_n)^T)^T$.
%
%
%
%
We assume that the network of nodes is an undirected network $ {\cal G} = ({\cal V},{\cal E}), $ where $ {\cal V} $ is the set of nodes and $ {\cal E} $ is the set of edges,  i.e., all pairs $\{i,j\} $ of nodes which can exchange information through a communication link.

 \noindent {\bf Assumption A2.}
 The network  $ {\cal G} = ({\cal V},{\cal E}) $ is connected, undirected and simple (no self-loops nor multiple links).

 Let us denote by $ O_i $ the set of nodes that are connected with node $ i $ through
 a direct link (neighborhood set),
   and let $\bar{O}_{i}=O_{i}\bigcup \{i\}.$  Associate with $ {\cal G} $ a symmetric, doubly stochastic $n \times n$ matrix $ W. $ The elements of $ W $ are all nonnegative and both rows and columns sum up to one. More precisely,  the following is assumed.

 \noindent{\bf Assumption A3.} The matrix $ W = W^T \in \mathbb{R}^{n \times n} $ is doubly stochastic,
 with elements $ w_{ij} $ such that
 $$ w_{ij} > 0 \mbox{ if } \{i,j\} \in {\cal E}, \;  w_{ij} = 0 \mbox{ if }  \{i,j\} \notin {\cal E},\,i \neq j, \mbox{ and }  w_{ii} = 1 - \sum_{j \in O_i} w_{ij} $$
 and there exist constants $\wmin$ and $\wmax$ such that for $i=1,\ldots,n$
 $$0 < \wmin \leq w_{ii} \leq  \wmax <1.$$

 Denote  by $ \lambda_1 \geq \ldots \geq \lambda_n $ the eigenvalues of $ W. $
 It can be shown that $ \lambda_1 = 1,$ and $|\lambda_i|<1$, $i=2,...,n$.

For future reference, define the
$n \times n$ matrix $J$ that has all
the entries equal~$1/n$.
 We refer to $J$
 as the ideal consensus matrix;
 see, e.g.,~\cite{MouraConsensus}. 
 Also, introduce the
$(nd) \times (nd)$ matrix
 $\mathcal{W}=W \otimes I$,
 where $\otimes$ denotes the Kronecker product
  and $I $ is the identity matrix from  $ \mathbb{R}^{d \times d}. $
   It can be seen that $d \times d$ block on the
   $(i,j)$-th position
   of the matrix $\mathcal{W}$
   equals to $w_{ij}\,I$.
   By properties of the Kronecker
   product, the eigenvalues
   of $\mathcal{W}$ are
   $\lambda_1,...,\lambda_n$,
   each one occurring with the multiplicity~$d$.

\subsection{Exact Distributed first order method}
Let us now briefly review the distributed first order method in~\cite{harnessing};
see also~\cite{angelia,angelia2}. 
These methods serve as a basis
for the development of the proposed
distributed spectral gradient method.
The method in~\cite{harnessing}
maintains over iterations
$k=0,1,...,$ at each node $i$,
the solution estimate $x_i^k \in {\mathbb R}^d$
and an auxiliary variable $u_i^k \in {\mathbb R}^d$.
Specifically, the update rule is
as follows
\begin{eqnarray}
\label{eqn-harnessing-node-persp-1}
x_i^{k+1} &=&
 \sum_{j \in \bar{O}_i} w_{ij}\,x_j^{(k)} - \alpha\,u_i^k\\
 \label{eqn-harnessing-node-persp-2}
u_i^{k+1} &=&
\sum_{j \in \bar{O}_i} w_{ij}\,u_j^{(k)}
+ \left( \nabla f_i(x_i^{k+1}) -  \nabla f_i(x_i^{k}) \right),\,\,k=0,1,...
\end{eqnarray}
Here, $\alpha>0$ is a constant step-size; the initialization
$x_i^0$, $i=1,...,n$,
is arbitrary, while
 $u_i^0 = \nabla f_i(x_i^0)$, $i=1,...,n$.
  Equation \eqref{eqn-harnessing-node-persp-1}
  shows that each node $i$,
  as with standard distributed gradient method~\cite{nedic_T-AC},
  makes  two-fold progress: 1)
  by weight-averaging
  its solution estimate with its' neighbors;
  and 2) by taking a step
  opposite to the estimated gradient direction.
  The standard distributed gradient method in \cite{nedic_T-AC}
  takes a negative step in the direction of
  $\nabla f_i(x_i^k)$, while the method in \cite{harnessing}
   makes a step in direction of $ u_i^k. $ This vector
   serves as a tracker
   of the network-wide
   gradient
   $\sum_{i=1}^n \nabla f_i(x_i^k)$.
    This modification
    in the update rule
    enables
    convergence
    to the
   exact solution under a constant step-size~\cite{harnessing}.

It is useful to represent method \eqref{eqn-harnessing-node-persp-1}--\eqref{eqn-harnessing-node-persp-2}
in vector format.
Let $x^k  \in {\mathbb R}^{nd}$,
 $u^k  \in {\mathbb R}^{nd}$,
 and recall function $F$ in \eqref{dobjective} and matrix $\mathcal{W} = W \otimes I$.
 Then, the method~\eqref{eqn-harnessing-node-persp-1}--\eqref{eqn-harnessing-node-persp-2}
 in the vector form becomes
\begin{eqnarray}
\label{eqn-harnessing-vector-1}
x^{k+1} &=&\mathcal{W}\,x^{(k)} - \alpha\,u^k\\
 \label{eqn-harnessing-vector-2}
u^{k+1} &=&
\mathcal{W}\,u^{(k)}
+ \left( \nabla F(x^{k+1}) -  \nabla F(x^{k}) \right),\,\,k=0,1,...,
\end{eqnarray}
with arbitrary $x^0$
 and $u^0 = \nabla F(x^0)$.

The method \eqref{eqn-harnessing-vector-1}--\eqref{eqn-harnessing-vector-2}
 allows for a primal-dual interpretation; see~\cite{dusan} and also~\cite{angelia2} 
 for a similar interpretation. 
 The primal-dual interpretation
 will be important
 for the development of the proposed
 distributed spectral gradient method.
 Namely, it is
 demonstrated in \cite{dusan} that \eqref{eqn-harnessing-vector-1}--\eqref{eqn-harnessing-vector-2}
 is equivalent to the following update rule
\begin{eqnarray}
\label{eqn-harnessing-PD-1}
x^{k+1} &=& {\cal W} x^k - \alpha (\nabla F(x^k) + u^k)\\
\label{eqn-harnessing-PD-2}
u^{k+1} &=& {\cal W} u^k + ({\cal W} - {\cal I}) \nabla F(x^k),
\end{eqnarray}
with variable $ u^0 = 0  \in \mathbb{R}^{d n}$
and arbitrary $x^0. $
It can be shown that,
under appropriately
chosen step-size $\alpha$,
the sequence $\{x^k\}$ converges to $x^*:=\mathbf{1} \otimes {y^*} =
(\,({y^*})^T,...,({y^*})^T\,)^T$,
and $u^k$ converges to
$-\nabla F(\mathbf{1}\otimes {y^*})
 = -
 (\nabla f_1({y^*})^T,...,\nabla f_n({y^*})^T)^T$. Here, $\mathbf{1} \in \mathbb{R}^{n} $ is
the vector with all components equal to one.

\subsection{Centralized spectral gradient method}
Let us briefly review the spectral gradient (SG) method in centralized optimization. Consider
the unconstrained minimization problem with  a generic objective function $ \phi: \mathbb{R}^d \to \mathbb{R} $ which is continuously differentiable. Let the initial
solution estimate be arbitrary $ x^0 \in {\mathbb R}^d. $
The SG method generates the sequence of iterates
$\{x^k\}$ as follows
\be \label{sgiteration}
x^{k+1} = x^{k} - \sigma_k^{-1} \, \nabla \phi(x^k), \;   k= 0,1,\ldots,
\ee
where the initial spectral coefficient $ \sigma_0 > 0 $ is arbitrary and $ \sigma_k $,
$k=1,2,...,$ is
 given by
\be
\label{sigmak}
\sigma_k = \mathcal{P}_{[\,{\sigma_{\min}},{\sigma_{\max}}\,]}
(\sigma_k^\prime),\; \;
\sigma_k^\prime=\frac{(s^{k-1})^T y^{k-1}}{(s^{k-1})^T s^{k-1}}.
\ee
Here,
$0<{\sigma_{\min}}< {\sigma_{\max}}<+\infty$ are given constants,
$ s^{k-1} = x^{k} - x^{k-1}, \; y^{k-1} = \nabla \phi(x^{k}) - \nabla \phi(x^{k-1}) $,
and $\mathcal{P}_{[a,b]}$ stands
for the projection of a scalar onto the interval $[a,b]$.
The projection onto the interval
$[{\sigma_{\min}},{\sigma_{\max}}]$
is the safeguarding that is necessary
for convergence.
The spectral coefficient  $ \sigma_k^\prime $ is derived as follows.
 Assume that the Hessian approximation in the form $ B_{k} =  \sigma_{k}  I. $ Then the approximate secant equation
\be \label{secant}
B_{k} s^{k-1} \approx y^{k-1}
\ee
can be solved in the least square sense. It is easy to show that least squares solution of (\ref{secant}) yields exactly (\ref{sigmak}).
For future reference, we briefly
review the result on the evolution of
error for the SG method stated in \cite{raydan1993}.
 Consider the special case of a strongly convex quadratic function $ \phi(x) = \frac{1}{2} x^T A x + b^T x $ for a symmetric positive definite matrix $ A$, and denote by $e^k:={x^*}- x^k$ the error at iteration~$k$,
 where $x^\star$ is the minimizer of $\phi$.
  Then, it can be shown that the error evolution can be expressed as~\cite{raydan1993}:
\be \label{sgerror}
e^{k+1} = (I - \sigma_k^{-1} A) e^k .
\ee
The above relation will play a key role in the intuitive explanation of the distributed spectral gradient method proposed in this paper.

\subsection{Small gain theorem}
Convergence analysis of the proposed method will be based upon the Small Gain Theorem, e.g.~\cite{smallgain}. 
This technique has been previously used and proved successful for the analysis 
 of exact distributed gradient methods in, e.g., \cite{angelia,angelia2}.  
 We briefly introduce the concept here, while more details are available in \cite{smallgain, angelia}.

 Denote by $ \mathbf{a}:=a^1,a^2,\ldots $ an infinite sequence of vectors, $ a^k \in \mathbb{R}^d, \; k=0,1,\ldots. $ For a fixed $ \delta \in (0,1) $, define
$$
\|\mathbf{a}\|^{\delta, K} = \max_{k=0,1,\ldots, K} \{\frac{1}{\delta^k}\|a^k\|\}
$$
$$
\|\mathbf{a}\|^{\delta} = \sup_{k \geq 0} \{\frac{1}{\delta^k}\|a^k\|\}.
$$
Obviously, for any $ K^\prime \geq K \geq 0 $ we have $ \|\mathbf{a}\|^{\delta,K} \leq  \|\mathbf{a}\|^{\delta,K^\prime} \leq  \|\mathbf{a}\|^{\delta}.  $ Also, if $ \|\mathbf{a}\|^{\delta} $ is finite for some $ \delta \in (0,1) $ than the sequence $ \mathbf{a} $ converges to zero R-linearly. 
 We present the Small Gain Theorem in a simplified form that involves only two sequences, as this will  will suffice for our considerations; for more general forms of the result see \cite{smallgain, angelia}.
\begin{te} \label{tesmallgain} \cite{smallgain, angelia}.
Consider two infinite sequences $ \mathbf{a} = a^0,a^1,\ldots, \; \mathbf{b} = b^0,b^1,\ldots, $ with $ a^k, b^k \in \mathbb{R}^d, \; k=0,1,\ldots.$ Suppose that for some $ \delta \in (0,1)$ and for all $ K = 0,1,\ldots, $ there holds
$$
\|\mathbf{a}\|^{\delta,K} \leq \gamma_1 \|\mathbf{b}\|^{\delta, K} + w_1
$$
$$
\|\mathbf{b}\|^{\delta,K} \leq \gamma_2 \|\mathbf{a}\|^{\delta, K} + w_2,
$$
where $ \gamma_1,  \gamma_2 \in [0,1). $ Then
$$
\|\mathbf{a}\|^{\delta} \leq \frac{1}{1-\gamma_1 \gamma_2}(w_1 \gamma_2 + w_2).
$$
Furthermore,
$ \lim_{k \to \infty} a^k =0 $ R-linearly.
\end{te}

%
%

Following, for example, the proof of Lemma~6 in \cite{dusan} (see also \cite{angelia}), it is easy to derive  the  result below.
\begin{lemma} \label{lemadusan3}
Consider three infinite sequence $ \mathbf{a} = a^0,a^1,\ldots, \; \mathbf{b} = b^0,b^1,\ldots, \; \mathbf{c} = c^0,c^1,\ldots $ with $ a^k, b^k, c^k \in \mathbb{R}^d, \; k=0,1,\ldots.$ Suppose that there holds
$$ \|a^{k+1}\| \leq c_1 \|a^k\| + c_2 \|b^k\|+c_3 \|c^k\|,  \; k=0,1,\ldots $$
where $ c_1,c_2, c_3 \geq  0. $ Then, for all $ K= 0,1,\ldots $ and $0 \leq c_1 <\delta<1,$
$$ \|\mathbf{a}\|^{\delta,K} \leq  \frac{c_2}{\delta-c_1} \|\mathbf{b}\|^{\delta,K} +\frac{c_3}{\delta-c_1} \|\mathbf{c}\|^{\delta,K}+ \frac{\delta}{\delta-c_1}\|a^0\|. $$
\end{lemma}

\section{Spectral gradient method for distributed optimization}

\subsection{The algorithm}
Let us now present the proposed Distributed Spectral
Gradient, DSG,  method. The method incorporates spectral-like step size policy into \eqref{eqn-harnessing-vector-1}--\eqref{eqn-harnessing-vector-2}.
 The step-sizes are locally computed and
 vary both across nodes and
 across iterations.
 As~\eqref{eqn-harnessing-vector-1}--\eqref{eqn-harnessing-vector-2},
 the DSG method maintains the sequence of
 solution estimates
 $x^k \in {\mathbb R}^{nd}$
 and an auxiliary sequence
 $u^k \in {\mathbb R}^{nd}$.
  Specifically, the update rule
 is as follows
\begin{eqnarray}
\label{eqn-DSG-vector-1}
x^{k+1} &=& {\cal W}\, x^k - \Sigma_k^{-1}\,u^k\\
\label{eqn-DSG-vector-2}
u^{k+1} &=& {\cal W} \, u^k + \left( \nabla F(x^{k+1}) - \nabla F(x^{k}) \right),\,\,k=0,1,...
\end{eqnarray}
The initial solution estimate~$x^0$ is arbitrary,
while~$s^{0} = \nabla F(x^0)$.
Here,
$$
\Sigma_k = diag\left(\sigma_1^k\,I, \ldots, \sigma_n^k\,I\right),
$$
is the $nd \times nd$
diagonal matrix that
collects inverse step-sizes
$\sigma_i^k$ at all nodes $i=1,...,n$.
The inverse step-sizes $\sigma_i^k$ are given by:
\begin{eqnarray}
 \label{dsigmak2}
\sigma_i^k &=&
\mathcal{P}_{[{\sigma_{\min}},{\sigma_{\max}}]}
\left\{\,\frac{(s_i^{k-1})^T y_i^{k-1}}{(s_i^{k-1})^T s_i^{k-1}}
+ \sum_{j \in \bar{O}_i}w_{ij}\,
\left(1 -  \frac{(s_i^{k-1})^T y_j^{k-1}}{(s_i^{k-1})^T s_i^{k-1}}  \right)\,\right\}\\
 s_i^{k-1} &=& x_i^k - x_i^{k-1} \nonumber \\
 y_i^{k-1} &=& \nabla f_i(x_i^{k}) - \nabla f_i(x_i^{k-1}), \nonumber
\end{eqnarray}
where $0<{\sigma_{\min}}<{\sigma_{\max}}<+\infty$
 are, as before, the safeguarding parameters.

Notice that the proposed step-size choice
does not incur an additional
communication overhead;
each node $i$ only needs to additionally store in its
memory $u_j^k$ for all its neighbors
$j \in O_i$.

In view of \eqref{eqn-harnessing-PD-1}--\eqref{eqn-harnessing-PD-2}, the method~\eqref{eqn-DSG-vector-1}--\eqref{eqn-DSG-vector-2}
 can be equivalently represented as follows
\begin{eqnarray}
\label{dsg1}
x^{k+1} &=& {\cal W} x^k - \Sigma_k^{-1} (\nabla F(x^k) + u^k)\\
\label{dsg2}
u^{k+1} &=& {\cal W} u^k + ({\cal W} - {\cal I}) \nabla F(x^k),\,\,k=0,1,...,
\end{eqnarray}
with variable $u^0 = 0 \in {\mathbb R}^{nd}. $

At the beginning of each iteration $ k+1,  $ a node $ i $ holds the current $ x_i^k, \nabla f_i(x_i^k), u_i^k, $ computes $ s_i^k = x_i^k - x_i^{k-1}, \; y_i^{k-1} = \nabla f_{i}(x_i^k) - \nabla f_{i}(x_i^{k-1}) $ and computes $ \sigma_i^k $ by (\ref{dsigmak2}). After that, it updates its' estimation of $ x_i $ through communication with all neighbouoring nodes $ j \in O_i $ as
$$ x_i^{k+1} = \sum_{j \in \bar{O}_i} w_{ij} x_j^k - (\sigma_i^k)^{-1} \left( \nabla f_i(x_i^k) + u_i^k \right) $$
$$ u_i^k = \sum_{j \in \bar{O}_i} w_{ij} u_j^k + \sum_{j \in \bar{O}_i} w_{ij} \nabla f_j(x_j^k) - \nabla f_i(x_i^k) .$$
Therefore, the iteration is fully distributed and each node interchanges messages only locally, with immediate neighbouors.

We next comment on the safeguarding parameters in \eqref{dsigmak2}. In practice, 
  the safeguarding upper bound ${\sigma}_{\max}$ can be
 set to a large number, e.g., ${\sigma}_{\max} = 10^{8}$;
 the safeguarding lower bound can be set to ${\sigma}_{\min} = \frac{L}{c}$,
 with $c \in [10,100]$. This in particular means
 that the proposed algorithm~\eqref{dsg1}--\eqref{dsg2}
    can take step-sizes~$\frac{1}{\sigma_i^k}$ that are much
    larger than the maximal allowed step-sizes with~\cite{harnessing}.
    In other words, as shown in Section~5 by simulations,
    ${\sigma}_{\min}$ can be chosen such
    that the method in \cite{harnessing} with
    step-size $\alpha = 1/{\sigma}_{\min}$ diverges,
    while the novel method~\eqref{dsg1}--\eqref{dsg2}
    with time-varying step sizes and the safeguard lower bound ${\sigma}_{\min}$
     (hence potentially taking step-size values close or equal to $1/{\sigma}_{\min}$) still converges.
%
%
\subsection{Step-size derivation}

We now provide a derivation and a justification
of the step-size choice~\eqref{dsigmak2}.
For notational simplicity, assume for the rest of this Subsection that $ d =1 $ and thus $ {\cal W} = W $.
 Let each $f_i$ be a strongly convex quadratic function, i.e.,
$$ f_i(x_i) = \frac{1}{2} h_i(x_i - b_i)^2, $$ and $ H = diag(h_1,\ldots,h_n), $
$h_i>0$, for all $i$.
Then, for the
primal error $e^k := x^k-x^*$
and the dual error
$\tilde{u}^k:= u^k+\nabla F(x^*)$, one can show that the following
recursion holds:
\be \label{derror}
 \begin{bmatrix}  e^{k+1} \\ \tilde{u}^{k+1} \end{bmatrix} = \begin{bmatrix} W - \Sigma_k^{-1} H & - \Sigma_k^{-1} \\ (W - I)H & W - J \end{bmatrix} \cdot \begin{bmatrix} e^k \\ \tilde{u}^k
 \end{bmatrix}
\ee
We now make a parallel and identification
between the error dynamics of the centralized SG method for
a strongly convex quadratic cost with leading matrix $A$ given in (\ref{sgerror})
 and the error dynamics of the proposed distributed method in~(\ref{derror}).
 Consider first the centralized SG method.
 The error dynamics
 matrix is given by $ I - \sigma_k^{-1}  A $,
 while the (new) spectral coefficient
 is sought to fit the secant equation with least mean square deviation:
  $\sigma_{k} (x^{k} - x^{k-1}) = A (x^{k} - x^{k-1})$.
  That is, the error dynamics
  matrix $I - \sigma_k^{-1}  A$ is made small
  by letting
  $\sigma_k\,I$
  be a scalar matrix approximation for
  matrix $A$, i.e.,
  solving
  $$ \min_{\sigma > 0}  \|\sigma s^k - y^k \|^2 = \min_{\sigma > 0} \|\sigma s^k - A s^k \|^2. $$
  Now,
  consider the error dynamics
  of the proposed distributed method in \eqref{derror},
  and specifically focus on the
  update for the primal error:
  \begin{eqnarray}
  e^{k+1} &=& \left( W - \Sigma_k^{-1} H \right)\,e^k + \Sigma_k^{-1}\,\widetilde{u}^k \nonumber \\
  &=&
  \label{eqn-DSG-update-primal-error}
  \left(\, I - \Sigma_k^{-1} \,\left[ \Sigma_k\,(I-W)+H \right] \,\right)\,e^k + \Sigma_k^{-1}\,\widetilde{u}^k.
  \end{eqnarray}
  Notice  that the second error equation in (\ref{derror}) does not depend on $\Sigma_k$.
Comparing \eqref{sgerror} with \eqref{eqn-DSG-update-primal-error},
we first see that both
the primal and the dual error play a role in~\eqref{eqn-DSG-update-primal-error}.
 The effect of the dual error
 $\widetilde{u}^k$ can be controlled by making
 $\Sigma_k^{-1}$ small enough.
 This motivates the safeguarding of
 $\Sigma_k$ from below by~${\sigma_{\min}}$.
 Regarding the effect of the primal error
  $e^k$,
  one can see that the it influences the error through the  matrix~$ I - \Sigma_k^{-1} \,\left[ \Sigma_k\,(I-W)+H \right]$.
 Analogously
   to the centralized SG case, this matrix can be made small by
   the following identification
   $$A \equiv \Sigma_k\,(I-W)+H, \mbox{ and }
   \sigma_k \equiv \Sigma_k.$$
   Therefore, we seek
    $\Sigma_{k+1}$
   as the least mean squares
   error fit to the following equation
   $$
   \Sigma_{k+1}\, \left(x^{k+1} - x^k\right) =  \left(\Sigma_k (I - W) + H\right)\,\left(x^{k+1} - x^k\right).
   $$
  For generic (non-quadratic) cost functions,
  this translates into the following:
  $$
   \Sigma_{k+1} \,\left(x^{k+1} - x^k\right) =
   (\,\Sigma_k (I - W)\,)\,\left(x^{k+1} - x^k\right) +
   \left(\,\nabla F(x^{k+1}) - \nabla F(x^k)\,\right).
  $$
The (intermediate) inverse step-size matrix $\Sigma^\prime_{k+1}$
is now obtained by minimizing
\[
\left\|  \Sigma_{k+1} \,\left(x^{k+1} - x^k\right) -
   (\,\Sigma_k (I - W)\,)\,\left(x^{k+1} - x^k\right) -
   \left(\,\nabla F(x^{k+1}) - \nabla F(x^k)\,\right) \right\|^2.
\]
Finally, to
ensure strictly positive step-sizes on the one hand,
and a bounded effect of the dual error on the other hand,
 $\Sigma^\prime_{k+1}$ is projected
 entry-wise onto the interval $[\,{\sigma_{\min}},\,{\sigma_{\max}}\,]$.

\section{Convergence results}

In this section we will prove that the proposed DSG method, (\ref{dsg1})-(\ref{dsg2}), converges to the solution of problem (\ref{objective}) provided that the   spectral coefficients $\sigma^k_i$ are uniformly bounded with properly choosen constants. For the sake of simplicity, we will restrict our attention to one dimensional case, i.e., $d=1$, while the general case is proved analogously. Hence, we have $ {\cal W} = W $ in this Section.

The following notation and relations are used. Recall that $x^*= \1 \otimes y^*$ where $y^*$ is the solution of (\ref{objective}). Define $\tx^k=x^k-\jc x^k$ and $\bx ^k=\1 ^T x^k /n.$  Then
$$
 \tx^k=x^k-\frac{1}{n} \1 \1 ^T x^k=x^k-\1 \otimes \bx^k.
$$
Also, for  $e^k=x^k-x^*, $
$$ \label{41b} (I-\jc) e^k=(I-\jc) x^k- (I-\jc) x^*=\tx^k-x^*+\1 \otimes y^*=\tx^k.
$$
Moreover, notice that $\jc^2=\jc$ and therefore $\jc (I-\jc)=0,$ which further implies  $\jc \tx^k=0$.
Now, for  $\tilde{W}=W-\jc$ we obtain
$$
\tilde{W} \tx^k=(W-\jc)(I-\jc) x^k=W \tx^k-\jc \tx^k=W \tx^k
$$
and
\be \label{41d} (I-\jc) W e^k=W(I-\jc)  e^k=W \tx^k=\tilde{W} \tx^k.
\ee
Define $\eb^k=\bx ^k-y^*. $ So, the following equalities hold
\be \label{42} e^k=x^k- \1 \otimes {\bx}^k + \1 \otimes {\bx}^k- \1 \otimes y^*=\tx^k+\1 \otimes \eb^k.
\ee
Given that $W$ is doubly stochastic, there follows $W x^*=x^*$,    $\1^T W =\1^T$ and $\1^T (W-I)=0$. So, multiplying (\ref{dsg2}) from the left with $\1^T,$ we obtain $\bu^{k+1}=\bu^k,$ where $\bu^k=\1^T u^k /n$. Since  $u^k=0$,  we conclude that
\be
\label{43}
\bu^k=0, \quad k=0,1,...
\ee
See Lemma~8 in \cite{dusan} that applies here as well, since the update (\ref{dsg2}) is a special case of update (16) in \cite{dusan}, with ${\cal B}=0$ defined therein. Moreover, define $\tu^k=\nabla F(x^*)+u^k.$ Using the fact that $\1^T \nabla F(x^*)=0$ we obtain
\be \label{44} \jc \tu^k=\frac{1}{n} \1 \otimes (\1 ^T u^k+\1^T \nabla F(x^*))=\jc u^k=\1 \otimes \bu^k=0
\ee

Now, Assumption A1 together with the Mean value theorem implies that for all $i=1,2,...,n$ and $k=1,2,...$,  there exists $\theta^k_i$ such that
$$\nabla f_i (x_i^k)-\nabla f_i (y^*)=\nabla^2 f_i (\theta^k_i) (x_i^k-y^*).$$
Therefore, there exists  a diagonal matrix $H_k$ such that
\be \label{45} \nabla F(x^k)-\nabla F(x^*)=H_k (x^k-x^*)=H_k e^k, \quad \mu I \preceq H_k \preceq L I.
\ee
  The R-linear convergence result for the DSG method is stated in the following theorem. 
  The Theorem corresponds to a worst case analysis
     that does not take into account
     the specific form of $\sigma_{i}^k$ in \eqref{dsigmak2}
     but only utilizes information on the safeguarding parameters
     ${\sigma}_{\min}$ and ${\sigma}_{\max}$. 
       Hence, the Theorem may be seen as an extension of 
   Theorem~2 in~\cite{angelia} that assumes node-varying but time-invariant step-sizes (here step-sizes are both node- and time-varying),  though we
   follow here a somewhat different proof path.

\begin{te} \label{glavna} Suppose that the assumptions A1-A3 hold. There exist $ 0 < \sigma_{\min} < \sigma_{\max} $ such that the sequence $\{x^k\}_{k \in \mathbb{N}}$ generated by DSG method converges R-linearly to the solution of problem (\ref{objective}).
\end{te}

{\bf Proof.} Let us first introduce the notation
$ \sigma_{\min}^{-1} = \dmax $ and $ \sigma_{\max}^{-1} = \dmin, \; \Delta = \dmax - \dmin. $ Choose $ \dmin, \dmax $ such that
\be \label{new1}
\frac{\dmax}{\dmin} < 1 + \frac{\mu}{L}
\ee
and
\be \label{new2}
0 < {\dmin} < {\dmax} < \frac{1-\lambda_2}{\mu + L}.
\ee
Define $ \delta = \delta(\dmin,\dmax) $ such that  $ \delta >0 $ and
\be \label{new4a}
 1 > \delta > 1 - \dmin \mu + \Delta L.
 \ee
As $  1 - \dmin \mu + \Delta L < 1 $ due to (\ref{new1}), $ \delta(\dmin, \dmax) $ is well defined. To simplify notation from now we will write $ \delta $ to denote $ \delta(\dmin,\dmax). $
Due to (\ref{new2}) we have
$$ 1 - \dmin \mu + \Delta L > \lambda_2 +\dmax L > \lambda_2 $$ and hence
\be \label{new4}
0 < \delta - (1 - \dmin \mu + \Delta L) < \delta - (\lambda_2 +\dmax L) < \delta - \lambda_2.
\ee
Notice further that $ \delta - (1 - \dmin \mu + \Delta L) $ is a decreasing function of $ \dmax$ and $ \Delta $ and therefore decreasing $  \dmax, \Delta $ if needed does not violate \ref{new4} if (\ref{new1}-\ref{new2}) are satisfied. In fact one  can take $ \dmax, \Delta $ arbitrary small  with the corresponding $ \dmin $ without violating (\ref{new1})-(\ref{new4}).

Denote $D_k=\si$ and $d_i^k=(\sigma_i^k)^{-1}$ and notice that $d_i^k\geq \dmin $. Subtracting $x^*$ from both sides of (\ref{dsg1}) and using the fact that $W x^*=x^*$ we obtain
$$e^{k+1}=W e^k-D_k(\nabla F(x^k)+u^k \pm \nabla F(x^*))=W e^k-D_k(\nabla F(x^k)-\nabla F(x^*))-D_k \tu^k.$$
From (\ref{45}) we obtain
\be \label{47} e^{k+1}=(W-D_k H_k) e^k-D_k \tu^k.
\ee
Now, adding $\nabla F(x^*)$ on both sides of (\ref{dsg2}) we obtain
\bea \label{48} \tu^{k+1} &=& W u^k+(W-I)\nabla F(x^k)+\nabla F(x^*) \pm W \nabla F(x^*) \nonumber \\
&=& W\tu^k+(W-I)(\nabla F(x^k)-\nabla F(x^*)).
\eea
Using (\ref{44}) and (\ref{45}) we get
\be \label{49} \tu^{k+1}=(W-\jc) \tu^k+(W-I)H_k e^k.
\ee
Taking the norm and using (\ref{42}),  we obtain
\be \label{49b} \|\tu^{k+1}\|\leq \lambda_2 \| \tu^k\| +(1-\lambda_n)L(\|\tx^k\|+\sqrt{n} |\eb^k| ).
\ee
 Lemma \ref{lemadusan3} with $c_1=\lambda_2,\; c_2=c_3=(1-\lambda_n)L$ yields
\be \label{49c} \|\tu\|^{\delta, K} \leq \frac{c_2}{\delta-c_1}(\|\tx\|^{\delta, K} + |\sqrt{n} \eb|^{\delta, K})+\frac{\delta}{\delta-c_1} \| \tu^0\|,
\ee
with $ \delta-c_1 > 0 $ due to (\ref{new4}). Define $\gamma_1:=c_2/(\delta-c_1)$.

 Multiplying both sides of (\ref{dsg1}) from the left with $\frac{1}{n} \1^T$ and using $\1^T W=\1^T$, (\ref{43}) and $\1^T \nabla F(x^*) =0$ we obtain
\bea \label{50} \bx^{k+1} &=& \bx^k-\frac{1}{n} \sum_{i=1}^{n} d_i^k \nabla f_i(x_i^k)-\frac{1}{n} \sum_{i=1}^{n} d_i^k u_i^k \pm \frac{1}{n} \sum_{i=1}^{n} \dmin \nabla f_i(x_i^k) + \frac{1}{n} \sum_{i=1}^{n} \dmin u_i^k \nonumber \\
&=& \bx^k -\frac{\dmin}{n} \sum_{i=1}^{n}  \nabla f_i(x_i^k)-\frac{1}{n} \sum_{i=1}^{n} (d_i^k-\dmin) \nabla f_i(x_i^k)-\frac{1}{n} \sum_{i=1}^{n} (d_i^k-\dmin)u_i^k \nonumber \\
&\pm &  \frac{1}{n} \sum_{i=1}^{n} (d_i^k-\dmin) \nabla f_i(y^*) \nonumber \\
&=& \bx^k -\frac{\dmin}{n} (\sum_{i=1}^{n} \nabla f_i(x_i^k)-\1^T \nabla F(x^*)) \nonumber \\
&-& \frac{1}{n} \sum_{i=1}^{n} (d_i^k-\dmin) (\nabla f_i(x_i^k)-\nabla f_i(y^*)) -\frac{1}{n} \sum_{i=1}^{n} (d_i^k-\dmin)\tu_i^k.
\eea
Moreover, using arguments similar to (\ref{45}), we conclude that there are $\hat{H}_i^k, \tilde{H}^k \in [\mu,L]$ such that
$$\nabla f_i(x_i^k)-\nabla f_i(y^*)=\nabla f_i(x_i^k)-\nabla f_i(y^*) \pm \nabla f_i(\bx^k)=\hat{H}_i^k \tx_i ^k+\tilde{H}^k \eb^k. $$
So, using the above equality  in the first sum in (\ref{50}), and the inequality  (\ref{45}) in the second sum,  with $H_i^k$ being the i-th diagonal component of $H_k$, after subtracting $y^*$ from both sides,  we obtain
\bea
\eb^{k+1} &=& \eb^k -\frac{\dmin}{n} \sum_{i=1}^{n}  (\hat{H}_i^k \tx_i ^k+\tilde{H}^k \eb^k)
- \frac{1}{n} \sum_{i=1}^{n} (d_i^k-\dmin) H_i^k e_i^k -\frac{1}{n} \sum_{i=1}^{n} (d_i^k-\dmin)\tu_i^k \nonumber \\
&=& \eb^k(1-\dmin \tilde{H}^k) -\frac{\dmin}{n} \sum_{i=1}^{n} \hat{H}_i^k \tx_i ^k \nonumber \\
&-& \frac{1}{n} \sum_{i=1}^{n} (d_i^k-\dmin) H_i^k e_i^k -\frac{1}{n} \sum_{i=1}^{n} (d_i^k-\dmin)\tu_i^k.
\eea
Since  (\ref{new1}) implies  $\dmin<1/L$, we have  $|1-\dmin \tilde{H}^k |\leq 1-\dmin \mu$. Moreover,
$$|\eb^{k+1}| \leq (1-\dmin \mu) |\eb^k|+\frac{\dmin}{n}L\|\tx^k\|_{1}+\frac{\ddif}{n}(L\|e^k\|_{1}+\|\tu^k\|_{1}).$$
Using the norm equivalence $\|\cdot\|_{1} \leq \sqrt{n} \|\cdot \|_{2}$, and multiplying both sides of the previous inequality with $\sqrt{n}$, we get
$$
\sqrt{n} |\eb^{k+1}| \leq (1-\dmin \mu) \sqrt{n} |\eb^k|+\dmin L\|\tx^k\| + \ddif (L\|e^k\| +\|\tu^k\| ).
$$
Furthermore, taking (\ref{42}) into account, the previous inequality becomes
\be \label{53b} \sqrt{n} |\eb^{k+1}| \leq (1-\dmin \mu+\ddif L) \sqrt{n} |\eb^k|+(\dmin+\ddif) L\|\tx^k\| +\ddif \|\tu^k\| .
\ee
Lemma \ref{lemadusan3} with $\tilde{c}_1=1-\dmin \mu+\ddif L,\; \tilde{c}_2=\dmax L,\; \tilde{c}_3=\ddif$ implies
\be \label{53c} |\sqrt{n} \eb|^{\delta, K} \leq \frac{1}{\delta-\tilde{c}_1}(\tilde{c}_2\|\tx\|^{\delta, K} +\tilde{c}_3 \|\tu\|^{\delta, K}+\delta |\sqrt{n} \eb^0|),
\ee
for $\delta \in (\tilde{c}_1,1)$. Notice that (\ref{new4}) implies that $\delta - \tilde{c}_1>0. $  Define
$$ \theta_2 = \frac{\tilde{c}_3}{\delta - \tilde{c}_1}, \;  \gamma_2 =\frac{\tilde{c}_2}{\delta - \tilde{c}_1}. $$
 Incorporating (\ref{49c}) into (\ref{53c}) and rearranging, we obtain
\be \label{53d} |\sqrt{n} \eb|^{\delta, K} \leq \frac{\gamma_2+\theta_2 \gamma_1}{1-\theta_2 \gamma_1} \|\tx\|^{\delta, K} +\frac{\theta_2 \delta \|\tu^0\|}{(\delta-c_1)(1-\theta_2 \gamma_1)} +\frac{ \delta |\sqrt{n} \eb^0|}{(\delta-\tilde{c}_1)(1-\theta_2 \gamma_1)} ,
\ee
provided that $\theta_2 \gamma_1<1.  $ This condition reads
\be
\label{new5}
 \frac{\Delta}{\delta - (1 - \dmin \mu + \Delta L)}\frac{(1-\lambda_n)L}{\delta - \lambda_2} < 1.
 \ee
 Clearly, there exists $ \delta, \dmin, \dmax $ such that for $ \dmax, \Delta $ small enough (\ref{new5}) holds as the left-hand side expression in (\ref{new5}) is increasing function of $ \dmax, \Delta $ and the corresponding $ \dmin $ satisfies (\ref{new2}).

Now, multiplying (\ref{47}) from the left with $I-\jc$ and using (\ref{41b}) and (\ref{41d}), we have
$$\tx^{k+1}=\tilde{W} \tx^k-(I-\jc) D_k H_ke^k-(I-\jc) D_k \tu^k.$$
Furthermore,  (\ref{42}) implies
$$\tx^{k+1}=(\tilde{W}-(I-\jc)D_k H_k ) \tx^k-(I-\jc) D_k H_k (\1 \otimes \eb^k)-(I-\jc) D_k \tu^k.$$
The inequality $\|\tilde{W}\| \leq \lambda_2$ yields
$$
\|\tx^{k+1}\|\leq(\lambda_2+\dmax L ) \|\tx^k\|+ \dmax L \sqrt{n} |\eb^k|+\dmax \|\tu^k\|.
$$
Again,  Lemma \ref{lemadusan3} with $\hat{c}_1=\lambda_2+\dmax L ,\; \hat{c}_2=\dmax L,\; \hat{c}_3=\dmax$,  implies
$$
\|\tx\|^{\delta,K}\leq \frac{\hat{c}_2}{\delta-\hat{c}_1} |\sqrt{n} \eb|^{\delta,K}+\frac{\hat{c}_3}{\delta-\hat{c}_1}\|\tu\|^{\delta,K}+\frac{\delta}{\delta-\hat{c}_1}\|\tx^0\|,
$$
with  $\delta - \hat{c}_1 >0$ due to (\ref{new4}).  Define $\gamma_3=\hat{c}_2/(\delta-\hat{c}_1)$ and $\theta_3=\hat{c}_3/(\delta-\hat{c}_1)$. Using (\ref{49c})  and rearranging, we obtain
\be \label{56} \|\tx\|^{\delta,K}\leq \frac{\gamma_3+\theta_3}{1-\theta_3 \gamma_1} |\sqrt{n} \eb|^{\delta,K}+
\frac{\theta_3 \delta }{(\delta-c_1) (1-\theta_3 \gamma_1)}\|\tu^0\|+\frac{\delta}{(\delta-\hat{c}_1)(1-\theta_3 \gamma_1)}\|\tx^0\|,
\ee
with $ \theta_3 \gamma_1 < 1 $ for $ \dmax $ small enough, due to the fact that
$$
\theta_3  \gamma_1 = \frac{\dmax}{\delta - (\lambda_2 + \dmax L)} \frac{(1-\lambda_n)L}{\delta - \lambda_2}
$$
is an increasing function of $ \dmax. $

Finally, considering (\ref{53d}),  (\ref{56}) and Theorem \ref{tesmallgain}, we conclude that $\tx^k$ and $\eb^k$ tend to zero R-linearly if
\be
\label{57}
 \frac{\gamma_2+\theta_2 \gamma_1}{1-\theta_2 \gamma_1} \frac{\gamma_3+\theta_3}{1-\theta_3 \gamma_1}<1.
\ee
The definition of $ \gamma_2 $ implies that it can be  arbitrary small if $  \dmax $ is small enough. As already stated, $ \theta_2 \gamma_1/(1-\theta_2 \gamma_1) $ is increasing function of $  \Delta $ Therefore, taking $ \Delta $ small enough, with the proper choice of $ \dmin, $  one can make the first term in (\ref{57}) arbitrary small.  On the other hand,
$$ \theta_3 + \gamma_3 = \frac{\dmax(L+1)}{\delta - (\lambda_2 + \dmax L)} $$ is again increasing function of $\dmax$ as is the function $ (1- \theta_3 \gamma_1)^{-1}. $ So, for $ \dmax, \Delta $ small enough and $ \dmin $ such that  (\ref{new1}-\ref{new2}) hold,  the inequality (\ref{56}) holds and the statement is proved.  $\Box$

Theorem~{4.1} 
poses conservative requirements on safeguarding. However, 
 extensive simulations on strongly convex quadratic
    and logistic losses indicate that
    method \eqref{dsg1}--\eqref{dsg2}
    always converges for
    $\frac{1}{{\sigma}_{\min}}< c/L$ and
    $\frac{1}{{\sigma}_{\min}}> 1/\theta$,
    where $c$ and $\theta$ can be taken at
    least as large as $c=100$, $\theta=10^8$.

\section{Numerical experiments}
This section provides
a numerical example
to illustrate the performance of
the proposed distributed spectral method.
 The example demonstrates
 a significant speedup
 gained through
 the proposed spectral-like step-size policy
 with respect to the
 counterpart constant step-size method in~\cite{harnessing}.

Consider the problem with strongly convex local quadratic costs;
that is, for each $i=1,...,n$, let
 $f_i:\,{\mathbb R}^d \rightarrow \mathbb R$,
 $f_i(x)=\frac{1}{2}(x-b_i)^T A_i (x-b_i)$,
  $d=10$,
 where $b_i \in {\mathbb R}^d$
  and $A_i \in {\mathbb R}^{d \times d}$
   is a symmetric positive definite matrix.
   The data pairs $A_i,b_i$ are
   generated at random, independently across nodes,
   as follows. Each $b_i$'s entry is
   generated mutually independently
   from the uniform distribution on~$[1,31]$.
   Each $B_i$ is generated as
   $B_i = Q_i \,D_i\,Q_i^T$; here, $Q_i$ is the
   matrix of orthonormal eigenvectors
   of~$\frac{1}{2}(\widehat{B}_i+\widehat{B}_i^T)$,
   and $\widehat{B}_i$ is a matrix
   with independent, identically distributed~(i.i.d.)
   standard Gaussian entries;
    $D_i$ is a diagonal matrix
   with the diagonal entries drawn
   in an i.i.d. fashion from the uniform
   distribution on~$[1,101]$.

The network is a $n=30$-node
instance of the random geometric graph model
with the communication radius~$r =\sqrt{\frac{\mathrm{ln}(n)}{n}}$,
and it is connected.
The weight matrix~$W$ is
set as follows:
 for~$\{i,j\} \in E$, $i\neq j$, $w_{ij}=\frac{1}{2(1+\max\{d_i,d_j\})}$,
 where $d_i$ is the node $i$'s degree;
 for~$\{i,j\} \notin E$, $i\neq j$, $w_{ij}=0$;
 and $w_{ii}=1-\sum_{j \neq i}w_{ij}$, for all~$i=1,...,n$.

The
proposed DSG method is compared with the method in~\cite{harnessing}.
This is a meaningful
comparison as the method in \cite{harnessing}
is a state-of-the-art
distributed first order method, and
the proposed method is based upon it.
The comparison thus
allows to assess the benefits
of incorporating spectral-like
step-sizes into
distributed first order methods.
As an error metric, the relative error averaged across nodes
\[
\frac{1}{n}\sum_{i=1}^n \frac{\|x_i-{y^*}\|}{\| {y^*} \|},\,\,y^{*} \neq 0.
\]
is used.

All parameters for both algorithms
are set in the same way,
except for
step-sizes.
With the method in \cite{harnessing},
the  step-size is
$\alpha=1/(3L),$
where $L = \max_{i=1,...,n}\mu_i$,
and $\mu_i$ is the maximal
eigenvalue of~$A_i$.
This step-size corresponds
to the  maximal possible
step-size for the method in~\cite{harnessing} as empirically evaluated in \cite{harnessing}.
With the DSG method,
at all nodes
the initial step-size
value is set to $1/(3L)$. The safeguard
parameters on the step-sizes
are set to $10^{-8}$
(lower threshold for safeguarding),
and $10 \times \frac{1}{3L}$
 (upper threshold for safeguarding).
 Hence,
 the step-sizes in DSG
are allowed to
reach up to $10$
times larger values than
  the maximal possible
  value with the method from~\cite{harnessing}.

\begin{figure}
      \centering
      \includegraphics[height=3.9 in,width=2.9 in, angle = -90]{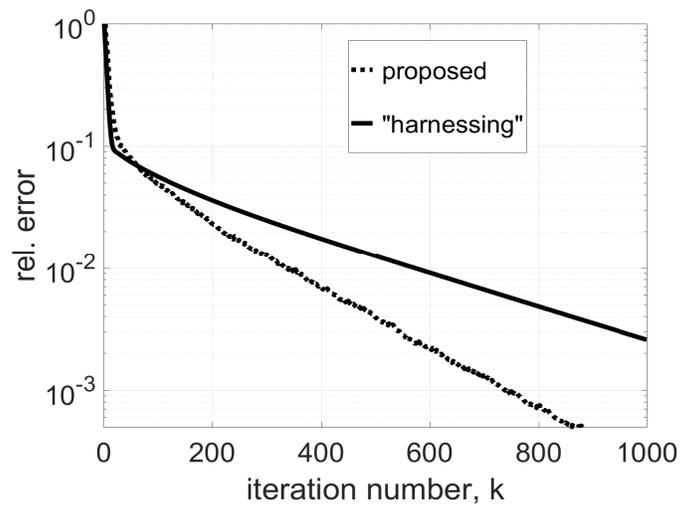}
      \includegraphics[height=3.9 in,width=2.9 in, angle = -90]{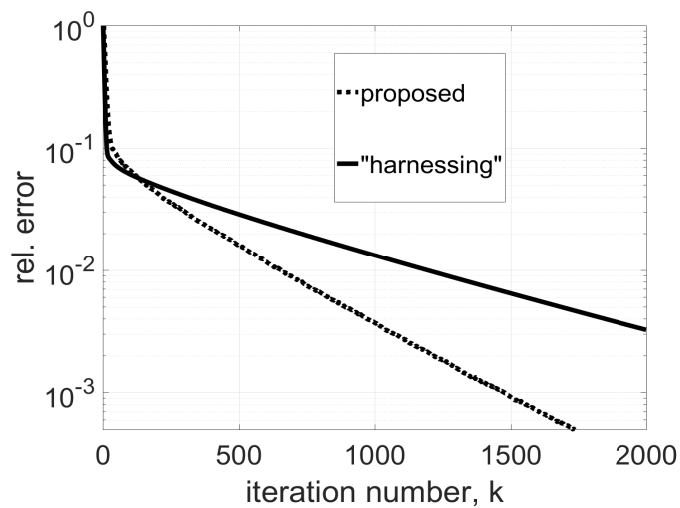}
      \caption{Relative error
      versus iteration number for the method in
      \cite{harnessing} (``harnessing'', solid line)
      and the proposed method (dotted line).}
      \label{FigTheoryExample}
\end{figure}

Figure~1 (top) plots
the relative error
versus number of iterations
with the two methods.One  can see that the DSG
method significantly
improves the convergence speed.
For example, to reach the relative error $0.01$,
the DSG method
requires about $340$
iterations, while
the method in \cite{harnessing}
takes about $560$
iterations for the same target accuracy;
this corresponds to
savings of about~$40\%. $

Figure~1 (bottom)
repeats the experiment
for a $n=100$-node
connected random geometric graph,
with the remaining data and network parameters
as before. We can see
that the DSG method achieves similar gains.
For example, for the $0.01$
accuracy, the DSG method needs about
$650$ iterations, while the method
in \cite{harnessing} needs about $1150$,
corresponding to decrease of about~$43\%$ in computational costs.
 We also report
that the method in \cite{harnessing} and step-size
equal to $1/{\sigma}_{\min}=10/(3L)$ diverges.

\section{Conclusion}

The method proposed in this paper, DSG, is a distributed version of the Spectral Gradient method for unconstrained optimization problems.   
Following the approach of exact distributed gradient methods in \cite{angelia} and \cite{harnessing}, 
at each iteration the nodes update two quantities -- the local approximation of the solution and the local approximation of the average gradient. The key novelty developed here is the step size selection which is defined in a spectral-like manner. Each node approximates the local Hessian by a scalar matrix thereby incorporating a degree of second order information in the gradient method. 
  The  spectral-like step-size coefficients are derived by  exploiting an analogy with the error dynamics of the classical spectral method for quadratic functions and embedding this dynamics into
 a primal-dual framework. This step-size calculation is computationally cheap and does not
 incur additional communication overhead. Under a set of standard assumptions regarding the objective functions and assuming a connected communication network, the DSG method generates a sequence of iterates which converges R-linearly to the exact solution of the aggregate objective function.
    The spectral gradient method is well known for its efficiency in classical, centralized optimization. Preliminary numerical tests demonstrate similar gains of incorporating spectral step-sizes in the distributed setting as well.


\begin{thebibliography}{99}
 
 \bibitem{bb} Barzilai J, Borwein JM, Two Point Step Size Gradient Methods, IMA Journal of
Numerical Analysis, 8 (1988), 141 - 148.

 
 \bibitem{bmr1} Birgin, E.G, Mart\'{\i}nez, J.M, Raydan M., Nonmonotone Spectral Projected Gradient
Methods on Convex Sets, SIAM Journal on Optimization, 10, (2000), 1196-1211.

\bibitem{bmr2} Birgin, E.G., Mart\'{\i}nez, J.M., Raydan M.,   Algorithm 813: SPG - Software for Convex-
Constrained Optimization,  ACM Transactions on Mathematical Software, 27 (2001), 340-349.

\bibitem{bmr3} Birgin, E.G., Mart\'{\i}nez, J.M., Raydan M Inexact Spectral Projected Gradient Methods
on Convex Sets, IMA Journal of Numerical Analysis, 23, (2003), 539-559.

\bibitem{bmr4} Birgin, E.G., Mart\'{\i}nez, J.M., Raydan M  Spectral Projected Gradient Methods: Review and Perspectives, Journal of Statistical Software 60(3), (2014), 1-21.

 \bibitem{BoydADMM}
Boyd, S., Parikh, N., Chu, E., Peleato, B., Eckstein, J., Distributed
  optimization and statistical learning via the alternating direction method of
  multipliers, Foundations and Trends in Machine Learning,
Volume 3, Issue 1, (2011) pp. 1-122.



 \bibitem{SayedEstimation}
Cattivelli, F.,  Sayed,  A.~H., Diffusion {LMS} strategies for distributed
  estimation, \emph{IEEE Transactions on  Signal  Processing},  vol.~58, no.~3, (2010) pp.
  1035--1048.



\bibitem{dl} Dai, Y.H., Liao, L.Z., R-Linear Convergence of the Barzilai and Borwein Gradient
Method,  IMA Journal on Numerical Analysis, 22 (2002),  1-10.

\bibitem{smallgain} Desoer, C., Vidyasagar, M., Feedback Systems: Input-Output Properties, SIAM 2009.

\bibitem{dusan} Jakoveti\'c, D., A Unification and Generalization of Exact Distributed First Order Methods,
 arxiv preprint, arXiv:1709.01317, 2017.


  \bibitem{arxivVersion}
Jakoveti\'c, D., Xavier, J.,  Moura, J.~M.~F.,  Fast distributed gradient
  methods, \emph{IEEE Transactions on  Automatic Control}, vol.~59, no.~5, (2014)  pp. 1131--1146.

\bibitem{cdc-submitted}
Jakoveti\'c, D.,  Moura, J.~M.~F., Xavier, J., Distributed {N}esterov-like
  gradient algorithms, in \emph{CDC'12, 51$^\textrm{st}$ IEEE Conference on
  Decision and Control}, Maui, Hawaii, December 2012, pp. 5459--5464.
  
  \bibitem{DQNSIAM} Jakoveti\'c, D., Bajovi\'c, D., Kreji\'c, N., Krklec Jerinki\'c, N., Newton-like Method with Diagonal Correction for Distributed Optimization,  SIAM J. Optimization, 27 (2), (2017), 1171-1203.



  \bibitem{SoummyaEst}
Kar, S.,  Moura, J.~M.~F., Ramanan, K., Distributed parameter estimation in
  sensor networks: Nonlinear observation models and imperfect communication,
  \emph{IEEE Transactions on Information Theory},  vol.~58, no.~6, (2012) pp.
  3575--3605.

  \bibitem{MouraConsensus}
Kar, S.,  Moura, J.~M.~F., Distributed Consensus Algorithms in Sensor Networks With Imperfect Communication: 
Link Failures and Channel Noise, 
\emph{IEEE Transactions on Signal Processing},  vol.~57, no.~1, (2009) pp.
  355--369.
  
  
\bibitem{SayedConf}
Lopes, C.,  {S}ayed,  A.~H., Adaptive estimation algorithms over distributed
  networks, in \emph{21st IEICE Signal Processing Symposium}, Kyoto, Japan,
  Nov. 2006.
  
  \bibitem{ribeiroNNpart1} Mokhtari, A., Ling, Q., Ribeiro, A., Network Newton--Part I: Algorithm and Convergence,
2015, available at: http://arxiv.org/abs/1504.06017


\bibitem{DQM} A. Mokhtari, W. Shi, Q. Ling, and A. Ribeiro,
DQM: Decentralized Quadratically Approximated Alternating Direction Method of Multipliers,
to appear in IEEE Trans. Sig. Process., 2016, {DOI}: 10.1109/TSP.2016.2548989

\bibitem{ESOM} A. Mokhtari, W. Shi, Q. Ling, and A. Ribeiro,
 A Decentralized Second Order Method with Exact Linear Convergence Rate for Consensus Optimization,
 2016, available at: http://arxiv.org/abs/1602.00596

 \bibitem{JoaoMotaMPC}
Mota, J., Xavier, J., Aguiar, P., P\"uschel, M., Distributed optimization
  with local domains: {A}pplications in {MPC} and network flows, \emph{to
  appear in IEEE Transactions on  Automatic Control}, 2015.

\bibitem{angelia} Nedic, A., Olshevsky, A., Shi, W., Uribe, C.A., Geometrically convergent distributed optimization with uncoordinated step-sizes, 
arXiv preprint, arXiv:1609.05877, 2016.

\bibitem{angelia2}  Nedic, A., Olshevsky, A., Shi, W., Achieving Geometric Convergence for Distributed Optimization over Time-Varying Graphs,
arxiv preprint, 	arXiv:1607.03218, 2016.


\bibitem{nedic_T-AC}
Nedi\'c, A., Ozdaglar, A., Distributed subgradient methods for multi-agent
  optimization, \emph{IEEE Transactions on Automatic Control},  vol.~54,
  no.~1, (2009) pp. 48--61.


\bibitem{harnessing} Qu, G., Li, N., Harnessing smoothness to accelerate distributed optimization, IEEE Transactions on Control of Network Systems (to appear)


 \bibitem{raydan1993} Raydan,  M.,  On the Barzilai and Borwein Choice of Steplength for the Gradient
Method,  IMA Journal of Numerical Analysis, 13 (1993),  321- 326.

\bibitem{raydan1997} Raydan M,  Barzilai and Borwein Gradient Method for the Large Scale Unconstrained
Minimization Problem, SIAM Journal on Optimization 7 (1997), 26 - 33.



\bibitem{RibeiroADMM1}
Schizas, I.~D. , Ribeiro, A.,  Giannakis,  G.~B., Consensus in ad hoc {WSN}s
  with noisy links -- {P}art~{I}: {D}istributed estimation of deterministic
  signals, \emph{IEEE Transactions on  Signal Processing},  vol.~56, no.~1, (2009) pp. 350--364.
  
  


\bibitem{extra}
Shi, W., Ling, Q., Wu, G., Yin, W.,  EXTRA: an Exact First-Order Algorithm for Decentralized Consensus Optimization,
\emph{SIAM Journal on Optimization}, No. 25 vol. 2, (2015) pp. 944-966.





 \end{thebibliography}
\end{document}